\documentclass[11pt]{article}
\usepackage{amsfonts}
\usepackage{mathrsfs}
\usepackage{amsmath}
\usepackage{amssymb}
\usepackage[mathscr]{eucal}
\usepackage{graphicx}
\renewcommand{\paragraph}{\roman{paragraph}}

\def \ov{\overline}

\newtheorem{theorem}{\scshape \mdseries  Theorem}[section]
\newtheorem{lemma}[theorem]{\scshape \mdseries  Lemma}

\begin{document}

\title{\sf Spectral Radius and Hamiltonicity of graphs}
\author{Guidong Yu$^{1}$\thanks{Email: guidongy@163.com.
Supported by the National Natural Science Foundation of China under
Grant no. 11371028, the Natural Science Foundation of Department of
Education of Anhui
Province of China under Grant nos. KJ2015ZD27, KJ2017A362.}, Yi Fang$^{1}$, Yizheng Fan$^{2}$, Gaixiang Cai$^{1}$\\
  {\small  \it $1.$ School of Mathematics \& Computation Sciences, Anqing Normal University, Anqing 246133, China}\\
{\small \it $2.$ School of Mathematical Sciences, Anhui University,
Hefei 230039, China }  }
\date{}
\maketitle

\noindent {\bf Abstract:} In this paper, we study the Hamiltonicity
of graphs with large minimum degree. Firstly, we present some
conditions for a simple graph to be Hamilton-connected and traceable
from every vertex in terms of the spectral radius of the graph or
its complement respectively. Secondly, we give the conditions for a
nearly balanced bipartite graph to be traceable in terms of spectral
radius, signless Laplacian spectral radius of the graph or its
quasi-complement respectively.


\noindent {\bf MR Subject Classifications:}  05C50, 05C45, 05C35

\section{Introduction}
Let $G=(V(G),E(G))$ be a simple graph of order $n$ with vertex set
$V(G)=\{v_1,v_2,\ldots,v_n\}$ and edge set $E(G)$. Denote by $e(G)=|
E(G)|$ the number of edges of the graph $G$. Let $N_G(v)$ be the set
of vertices which are adjacent to $v$ in $G$. The degree of $v$ is
denoted by $d_{G}(v)=|N_G(v)|$ (or simply $d(v)$), the minimum
degree of $G$ is denoted by $\delta(G)$. Let $X\subseteq V(G)$,
$G-X$ is the graph obtained from $G$ by deleting all vertices in
$X$. $G$ is called {\it $k$-connected} (for $k\in \mathbb{N}$) if
$|V(G)|>k$ and $G-X$ is connected for every set $X\subseteq V(G)$
with $|X| < k$. We note that $G$ is $k$-connected when $\delta
(G)\geq k$. A {\it regular graph} is one graph whose vertices all
have the same degrees, and a {\it bipartite semi-regular graph} is a
bipartite graph for which the vertices in the same part have the
same degrees. The {\it complement} of $G$ is denoted by
$\ov{G}=(V(\ov{G}),E(\ov{G}))$, where $V(\ov{G})=V(G)$,
$E(\ov{G})=\{xy:~x,y\in V(G), xy\not \in E(G)\}$. Let $G =(X, Y; E)$
be a bipartite graph with two part sets $X,Y$. If $|X|=|Y|$, $G =(X,
Y; E)$ is called a {\it balanced bipartite graph}. If $|X|=|Y|-1$,
$G =(X, Y; E)$ is called a {\it nearly balanced bipartite graph}.
The {\it quasi-complement} of $G =(X, Y; E)$ is denoted by
$\widehat{G} :=(X, Y; E')$, where $E' =\{xy : x\in X, y\in Y, xy\not
\in E \}$. For two disjoint graphs $G_1$ and $G_2$, the union of
$G_1$ and $G_2$, denoted by $G_1+G_2$, is defined as
$V(G_1+G_2)=V(G_1)\cup V(G_2)$ and $E(G_1+G_2)=E(G_1)\cup E(G_2)$;
and the join of $G_1$ and $G_2$, denoted by $G_1\vee G_2$, is
defined as $V(G_1\vee G_2)=V(G_1)\cup V(G_2)$, and $E(G_1\vee
G_2)=E(G_1+G_2)\cup\{xy:x\in V(G_1),y\in V(G_2)\}$. Denote $K_{n}$
the complete graph on $n$ vertices, $O_{n}=\ov K_{n}$ the empty
graph on $n$ vertices (without edges), $K_{n,m}=O_{n}\vee O_{m}$ the
complete bipartite graph with two parts having $n,m$ vertices, $G-v$
($v\in V(G)$) the graph obtained from $G$ by deleting $v$,
respectively.

The {\it adjacency matrix} of $G$ is defined to be a matrix
$A(G)=[a_{ij}]$ of order $n$,
  where $a_{ij}=1$ if $v_{i}$ is adjacent to $v_{j}$, and $a_{ij}=0$ otherwise.
The degree matrix of $G$ is denoted by
$D(G)=\hbox{diag}\left(d_G(v_1),d_G(v_2), \ldots, d_G(v_n)\right)$.
The matrix $Q(G)=D(G)+A(G)$ is the {\it signless Laplacian matrix}
(or {\it $Q$-matrix}) of $G$. Obviously, $A(G)$ and $Q(G)$  are real
symmetric matrix. So their eigenvalues are real number and can be
ordered. The largest eigenvalue of $A(G)$, denoted by $\mu(G)$, and
the corresponding eigenvectors (whose all components are positive
number) are called the {\it spectral radius} and the {\it Perron
vector} of $G$, respectively. The largest eigenvalue of $Q(G)$,
denoted by $q(G)$, is called the {\it signless Laplacian spectral
radius} of $G$.

A {\it Hamiltonian cycle} of the graph $G$ is a cycle of order $n$
contained in $G$, and a {\it Hamiltonian path} of $G$ is a path of
order $n$ contained in $G$, where $|V(G)|=n$. The graph $G$ is said
to be {\it Hamiltonian} if it contains a Hamiltonian cycle, and is
said to be {\it traceable} if it contains a Hamiltonian path. If
every two vertices of $G$ are connected by a Hamiltonian path, it is
said to be {\it Hamilton-connected}. A graph $G$ is {\it traceable
from a vertex} $x$ if it has a Hamiltonian $x$-path. The problem of
deciding whether a graph is Hamiltonian is one of
   the most difficult classical problems in graph theory.
Indeed, determining whether a graph is Hamiltonian is NP-complete.

Recently, the spectral theory of graphs has been applied to this
problem. Up to now, there are some references on the spectral
conditions for a graph to be traceable, Hamiltonian,
Hamilton-connected or traceable from every vertex. We refer readers
to see \cite{ning1, fied, yu1,  Mei, Li, Li1,  Liu, Ning0, Ning,
nikiforov, yu, yu4, yu2, yu3, zhou}. Particularly, Li and Ning
\cite{ning1} and Nikiforov \cite{nikiforov} study spectral
sufficient conditions of graphs with large minimum degree. Li and
Ning \cite{ning1} present some (signless Laplacian) spectral radius
conditions for a simple graph and a balanced bipartite graph to be
traceable and Hamiltonian, respectively. Nikiforov \cite{nikiforov}
gives some spectral radius conditions for a simple graph to be
traceable and Hamiltonian, respectively. Motivated by those papers,
in this paper, we also study the graphs with large minimum degree.
We will respectively present some conditions for a simple graph to
be Hamilton-connected and traceable from every vertex in terms of
the spectral radius of the graph or its complement in section 2, and
respectively give the conditions for a nearly balanced bipartite
graph to be traceable in terms of spectral radius, signless
Laplacian spectral radius of the graph or its quasi-complement in
section 3 .

\section{Spectral radius conditions for a graph to be Hamilton-connected,
and traceable from every vertex}

For an integer $k\geq0$, the $k$-closure of a graph $G$, denoted by
$C_k(G)$, is the graph obtained from $G$ by successively joining
pairs of nonadjacent vertices whose degree sum is at least $k$ until
no such pair remains, see [2]. The $k$-closure of the graph $G$ is
unique, independent of the order in which edges are added. Note that
$d_{C_k(G)}(u)+d_{C_k(G)}(v)\leq k-1$ for any pair of nonadjacent
vertices $u$ and $v$ of $C_k(G)$.

\begin{lemma} {\rm(Ore \cite{ore}, Bondy and Chv\'{a}tal \cite{bondy})}
{\rm (i)} If $G$ is a 2-connected graph of order $n$ and
$d_G(u)+d_G(v)\geq n+1$ for any two distant nonadjacent vertices $u$
and $v$, then $G$ is Hamilton-connected .

 {\rm (ii)} A 2-connected graph $G$ is Hamilton-connected if and only if $C_{n+1}(G)$ is so.
\end{lemma}

\begin{lemma} {\rm(Yu, Ye and Cai \cite{yu4})}
Let $G$ be a simple graph, with degree sequence
$(d_G(v_1),d_G(v_2),\\ \ldots,d_G(v_n))$, where $d_G(v_1)\leq
d_G(v_2)\leq\ldots\leq d_G(v_n)$ and $n\geq 3$. Suppose that there
is no integer $2\leq k\leq\frac{n}{2}$ such that $d_G(v_{k-1})\leq
k$, and $d_G(v_{n-k})\leq n-k$, then $G$ is Hamilton -connected.
\end{lemma}

\begin{lemma} {\rm(Hong and Shu \cite{hong1}, Nikiforov \cite
{nikiforov1})}
 If $G$ is a graph of order $n$, with $m$ edges and
minimum degree $\delta$, then
$$\mu(G)\leq\frac{\delta-1}{2}+\sqrt{2m-n\delta+\frac{(\delta+1)^2}{4}}.$$
\end{lemma}

\begin{lemma} {\rm(Hong and Shu \cite{hong1}, Nikiforov \cite
{nikiforov1})} If $2m\leq n(n-1)$, the function
$$f(x)=\frac{x-1}{2}+\sqrt{2m-nx+\frac{(x+1)^2}{4}}$$ is
decreasing in $x$ for $x\leq n-1$.
\end{lemma}

\begin{lemma} {\rm(Bondy and Murty \cite{bondy1})} Let $G$ be a graph. Then $G$ is traceable from
every vertex if and only if $G\vee K_1$ is Hamilton-connected.
\end{lemma}

Given a graph $G$ of order $n$, a vector $\textbf{x}\in R^{n}$ is
called to be defined on $G$, if there is a 1-1 map $\varphi$ from
$V(G)$ to the entries of $\textbf{x}$; simply written
$x_{u}=\varphi(u)$ for each $u\in V(G)$. If $\textbf{x}$ is an
eigenvector of $A(G)$, then $\textbf{x}$ is defined on $G$
naturally, $x_{u}$ is the entry of $\textbf{x}$ corresponding to the
vertex $u$. One can find that
$$\textbf{x}^TA(G)\textbf{x}=2\sum_{uv\in E(G)}x_ux_v,\eqno(2.1)$$ when $\mu$ is a
eigenvalue of $G$ corresponding to the eigenvector $\textbf{x}$ if
and only if $\textbf{x}\neq \mathbf{0}$,
$$\mu x_{v}=\sum_{u\in N_{G}(v)}x_{u},  \eqno(2.2)$$ for each vertex $v\in V(G)$. Equation (2.2) is called the eigenvalue-equation for the
graph $G$. In addition, for an arbitrary unit vector $\textbf{x}\in
R^n$,$$\mu(G)\geq \textbf{x}^TA(G)\textbf{x}, \eqno(2.3)$$ with
equality holds if and only if $\textbf{x}$ is an eigenvector of
$A(G)$ according to $\mu(G)$.

\begin{lemma} {\rm (Li and Ning \cite{ning1})}
Let $G$ be a graph with non-empty edge set. Then $$\mu(G)\geq
min\{\sqrt{d(u)d(v)}:uv\in E(G)\}. \eqno(2.4)$$ Moreover, if $G$ is
connected, then equality holds if and only if $G$ is regular or
bipartite semi-regular graph.
\end{lemma}

\begin{lemma} Let $G$ be a graph of order $n$. Then $$\mu(G\vee
K_1)>\frac{n-1}{n}\mu(G)+2\frac{\sqrt{n-1}}{n}.$$
\end{lemma}

{\bf Proof.} Let $\textbf{x}\in R^n$ be a unit Perron vector of $G$,
then by (2.1) and (2.3),
$$\mu(G)=\textbf{x}^TA(G)\textbf{x}=2\sum_{uv\in E(G)}x_ux_v.$$ Let
$w\in V(K_1), H=G\vee K_1,$ and let $\textbf{x}^{'}\in R^{n+1},
x_u^{'}=\sqrt{\frac{n-1}{n}}x_u$, for every $u\in V(G)$,
$x_w^{'}=\frac{1}{\sqrt{n}}.$ Since $\sum\limits_{u\in
V(G)}x_{u}^{2}=1,x_u> 0$, we have
$$\sum_{u\in V(H)}{x_u'}^2=\sum_{u\in V(G)}{x_u'}^2+{x_w'}^2=\frac{n-1}{n}\sum_{u\in V(G)}x_u^2+\frac{1}{n}=1,$$
and $\sum\limits_{u\in V(G)}x_u>\sum\limits_{u\in V(G)}x_u^2=1$.
Then by (2.1) and (2.3)
\begin{eqnarray*}
\mu(G\vee K_1)=\mu(H)&\geq& {\textbf{x}^{'}}^TA(H)\textbf{x}^{'}\\
&=& 2\sum_{uv\in E(G)}x_u^{'}x_v^{'}+2x_w^{'}\sum_{u\in V(G)}x_u^{'}\\
&=&2\frac{n-1}{n}\sum_{uv\in
E(G)}x_ux_v+2\frac{1}{\sqrt{n}}\sqrt{\frac{n-1}{n}}\sum_{u\in
V(G)}x_u\\
&>& \frac{n-1}{n}\mu(G)+2\frac{\sqrt{n-1}}{n}.
\end{eqnarray*}

So the result follows. \hfill $\blacksquare$

\begin{lemma} {\rm(Tomescu \cite{tomescu})}
Every $t$-regular graph on $2t~(t\geq3)$ not isomorphic to
$K_{t,t}$, or of order $2t+1$ for even $t\geq4$, is
Hamilton-connected.
\end{lemma}

\begin{lemma} Let $k\geq2,n\geq2k^2+1$, and $G$ be a graph of
order $n$. If $G$ is a subgraph of $K_2\vee(K_{n-k-1}+K_{k-1})$,
with minimum degree $\delta(G)\geq k$. Then $\mu(G)<n-k$, unless
$G=K_2\vee(K_{n-k-1}+K_{k-1})$.
\end{lemma}

{\bf Proof.} Set for short $\mu := \mu(G)$, and let
$\textbf{x}=(x_{v_1},\ldots,x_{v_n})^{T}$ be a unit Perron vector of
$G$. By (2.3), we have that
$$\mu=\textbf{x}^TA(G)\textbf{x}.$$

Assume that $G$ is a proper subgraph of
$K_2\vee(K_{n-k-1}+K_{k-1})$. By Perron-Frobenius theorem, we can
assume that $G$ is obtained by omitting just one edge $uv$ of
$K_2\vee(K_{n-k-1}+K_{k-1})$.

Write $X$ for the set of vertices of $K_2\vee(K_{n-k-1}+K_{k-1})$ of
degree $k$, let $Y$ be the set of their neighbors not in the set
$X$, and let $Z$ be the set of the remaining $n-k-1$ vertices of
$K_2\vee(K_{n-k-1}+K_{k-1})$.

Since $\delta(G)\geq k$, we can see that $G$ must contain all the
edges between $X$ and $Y$. Therefore, $\{u,v\}\subset Y\cup Z$, with
three possible cases: $(a)~\{u,v\}\subset Y; (b)~u\in Y, v\in Z;
(c)~ \{u,v\}\subset Z$. We shall show that case $(c)$ yields a graph
of no smaller spectral radius than case $(b)$, and that case $(b)$
yields a graph of no smaller spectral radius than case $(a)$.

Indeed, by (2.2), we have $x_i=x_j$ for any $i,j\in X$; likewise,
$x_i=x_j$ for any $i,j\in Y\backslash\{u,v\}$ and for any $i,j\in
Z\backslash\{u,v\}$. Thus, let
\begin{eqnarray*}
x &:=&x_i,i\in X,\\
y &:=&x_i,i\in Y\backslash\{u,v\},\\
z &:=&x_i,i\in Z\backslash\{u,v\}.
\end{eqnarray*}

Suppose that case $(a)$ holds, that is, $\{u,v\}\subset Y$. Choose a
vertex $w\in Z$, remove the edge $vw$ and add the edge $uv$.  Then
the obtained graph $G'$ is covered by case $(b)$.

If $x_w\leq x_u$, we have
$$\textbf{x}^TA(G')\textbf{x}-\textbf{x}^TA(G)\textbf{x}=2x_v(x_u-x_w)\geq0;$$

If $x_w>x_u$, swap the entries $x_u$ and $x_w$, write $\textbf{x}'$
for the resulting vector. We note that $\textbf{x}'$ is also a unit
vector, and have that
$$\textbf{x}'^TA(G')\textbf{x}'-\textbf{x}^TA(G)\textbf{x}=2(x_w-x_u)\sum_{i\in X}x_i\geq0.$$ Then by (2.3),
$\mu(G')\geq\mu(G)$, as claimed.

Essentially the same argument proves that case $(c)$ yields a graph
of no smaller spectral radius than case $(b)$. Therefore, we may
assume that $\{u,v\}\subset Z$. Since the vertices $u$ and $v$ are
symmetric, so $x_u=x_v$. Set $t :=x_u$ and note that the $n$
eigenvalue-equations of $G$ are reduced to four equations involving
just the unknowns $x,y,z,$ and $t$:
\begin{eqnarray*}
\mu x &=& (k-2)x+2y,\\
\mu y &=& (k-1)x+y+(n-k-3)z+2t,\\
\mu z &=& 2y+(n-k-4)z+2t,\\
\mu t &=& 2y+(n-k-3)z.
\end{eqnarray*} We find that
\begin{eqnarray*}
x&=&\frac{2y}{\mu-k+2},\\
z&=&(1-\frac{2(k-1)}{(\mu+1)(\mu-k+2)})y,\\
t&=&\frac{\mu+1}{\mu+2}(1-\frac{2(k-1)}{(\mu+1)(\mu-k+2)})y.
\end{eqnarray*}

Furtherly, note that if we delete all edges incident to vertices in
$X$, and add the edge $uv$ to $G$, we obtain the graph
$K_{n-k+1}+\ov K_{k-1}$. Letting $\textbf{x}''$ be the restriction
of $\textbf{x}$ to $K_{n-k+1}$, we find that
$$\textbf{x}''^TA(K_{n-k+1})\textbf{x}''=\textbf{x}^TA(G)\textbf{x}+2t^2-4(k-1)xy-(k-1)(k-2)x^2=\mu+2t^2-4(k-1)xy-(k-1)(k-2)x^2.$$

But since $\|\textbf{x}''\|^2=1-(k-1)x^2$, we see that
$$\mu+2t^2-4(k-1)xy-(k-1)(k-2)x^2=\textbf{x}''^TA(K_{n-k+1})\textbf{x}''\leq\mu(K_{n-k+1})\|\textbf{x}''\|^2=(n-k)(1-(k-1)x^2).$$

Assume for a contradiction that $\mu\geq n-k$. This assumption,
together with above inequality, yields $$\mu
+2t^2-4(k-1)xy-(k-1)(k-2)x^2\leq\mu(1-(k-1)x^2),$$ and therefore
$$2(k-1)xy-\frac{(\mu-k+2)(k-1)x^2}{2}\geq t^2.$$ Now, first combining above
equality about $x$, then combining about equality about $t$, we have
$$\frac{2(k-1)y^2}{\mu-k+2}\geq(\frac{\mu+1}{\mu+2})^2(1-\frac{2(k-1)}{(\mu+1)(\mu-k+2)})^2y^2.$$
Cancelling $y^2$ and applying Bernoulli's inequality to the right
side, we get
\begin{eqnarray*}
2(k-1)&\geq&(\mu-k+2)(1-\frac{1}{\mu+2})^2(1-\frac{2(k-1)}{(\mu+1)(\mu-k+2)})^2\\
&>&(\mu-k+2)(1-\frac{2}{\mu+2}-\frac{4(k-1)}{(\mu+1)(\mu-k+2)})\\
&=&\mu-k+2-\frac{2\mu-2k+4}{\mu+2}-\frac{4(k-1)}{\mu+1}\\
&>&\mu-k+2-\frac{2\mu+2k}{\mu+1}.
\end{eqnarray*}
Using the inequalities $\mu\geq n-k\geq2k^2-k+1$, we easily find
that $$2<\frac{2\mu+2k}{\mu+1}<3,$$ and so,
$$2(k-1)>2k^2-k+1-k+2-3=2k^2-2k,$$ a contradiction,
completing the proof. \hfill $\blacksquare$

\begin{theorem}
 Let $k\geq2, n\geq2k^2+1$ and let $G$ be a graph of order $n$ with
 minimum degree $\delta(G)\geq k$. If $$\mu(G)\geq n-k,$$ then $G$ is
 Hamilton-connected, unless $G=K_2\vee(K_{n-k-1}+K_{k-1})$.
\end{theorem}

{\bf Proof.} Assume that $\mu(G)\geq n-k$, but $G$ is not
Hamilton-connected. Let $H=C_{n+1}(G)$, then $H$ is not
Hamilton-connected by Lemma 2.1, $\delta(H)\geq\delta(G)\geq k,~{\rm
and}~ \mu(H)\geq\mu(G)\geq n-k$ by Perron-Frobenius theorem. Note
that $H$ is $(n+1)$-closure of $G$, thus every two nonadjacent
vertices $u,v$ have degree sum at most $n$, i.e.,$$d_H(u)+d_H(v)\leq
n. \eqno(2.5)$$ Since $H$ is not Hamilton-connected, by Lemma 2.2,
there is an integer $2\leq s\leq\frac{n}{2}$ such that
$d_H(v_{s-1})\leq s$ and $d_H(v_{n-s})\leq n-s$, obviously,
$s\geq\delta(H)\geq k$. Write $m$ for the number of edges of $H$,
set $\delta(H):=\delta$, then we can get

$$\begin{array}{rcl}
2m&=&\sum\limits_{i=1}^{s-1}d_H(v_i)+\sum\limits_{i=s}^{n-s}d_H(v_i)+\sum\limits_{i=n-s+1}^{n}d_H(v_i) \\
&\leq& s(s-1)+(n-2s+1)(n-s)+s(n-1)\\
&=&3s^2+n^2-2ns+n-3s.
\end{array} \eqno(2.6)$$

On the other hand, combining Lemmas 2.3, 2.4, we have
$$n-k\leq\mu(H)\leq\frac{k-1}{2}+\sqrt{2m-nk+\frac{(k+1)^2}{4}},$$
which, after some algebra operations, gives $$2m\geq
n^2-2kn+2k^2+n-2k. \eqno(2.7)$$ Next, we will prove that $s=k$.
Suppose $k+1\leq s\leq \frac{n}{2}$. Let $f(x)=3x^2+n^2-2nx+n-3x$,
we note $f(x)$ is convex in $x$, then $f(s)\leq f(k+1)$ or $f(s)\leq
f(\frac{n}{2}).$

Combining (2.6) and (2.7), we get $$n^2-2kn+2k^2+n-2k\leq2m\leq
f(s)\leq f(k+1)=3(k+1)^2+n^2-2n(k+1)+n-3(k+1)$$ or
$$n^2-2kn+2k^2+n-2k\leq 2m\leq f(s)\leq
f(\frac{n}{2})=\frac{3}{4}n^2-\frac{n}{2}.$$ Then
$n\leq\frac{k^2+5k}{2}$ or $n^2+(6-8k)n+8k(k-1)\leq0$, each of these
inequalities leads to a contradiction. So we have $s=k$, and thus
$\delta(H)=k$, then,
$$d_H(v_1)=d_H(v_2)=\ldots=d_H(v_{k-1})=k.$$

Our next goal is to show that $d_H(v_k)\geq n-k^2$. Indeed, suppose
that $$d_H(v_k)<n-k^2.$$ Also using Lemma 2.2, we get
\begin{eqnarray*}
2m&=&\sum_{i=1}^{k-1}d_H(v_i)+d_H(v_k)+\sum_{i=k+1}^{n-k}d_H(v_i)+\sum_{i=n-k+1}^{n}d_H(v_i)\\
&<&(k-1)k+n-k^2+(n-2k)(n-k)+k(n-1)\\
&=&n^2-2kn+2k^2+n-2k,
\end{eqnarray*}
contradicting (2.7). Hence $d_H(v_i)\geq n-k^2$ for every $i\in
\{k,k+1,\ldots,n\}.$

Next, we shall show that the vertices $v_{k},v_{k+1},\ldots,v_{n}$
induce a complete graph in $H$. Indeed, let
$v_i,v_j\in\{v_k,v_{k+1},\ldots,v_n\}$ be two distinct vertices of
$H$. If they are nonadjacent, then
\begin{eqnarray*}
\ d_H(v_i)+d_H(v_j)&\geq&2n-2k^2\\
&\geq& n+2k^2+1-2k^2\\
&=&n+1,
\end{eqnarray*}
contradicting (2.5).

Write $X$ for the vertex set $\{v_1,v_2,\ldots, v_{k-1}\}$. Write
$Y$ for the set of vertices in $\{v_k,v_{k+1},\ldots,v_n\}$ having
neighbors in $X$. Let $Z$ be the set of remaining vertices of
$V(G)$.

Since $|X|=k-1$, and $d_H(v_1)=d_H(v_2)=\ldots=d_H(v_{k-1})=k$, we
get $Y\neq\emptyset$, and any vertex in $X$ must have at least two
neighbors in $\{v_k,v_{k+1},\ldots,v_n\}$.

In fact, every vertex from $Y$ is adjacent to every vertex in $X$.
Indeed, suppose that this is not the case, and let
$w\in\{v_k,v_{k+1},\ldots,v_n\},u\in X,v\in X$, such that $w$ is
adjacent to $u$, but not to $v$. We see that $$d_H(w)+d_H(v)\geq
n-k+1+k=n+1,$$ contradicting (2.5).

Next, let $l=|Y|$ and note that $2\leq l\leq k$.

If $l=2$, then $H=K_2\vee(K_{n-k-1}+K_{k-1})$. Since $G\subseteq H$,
by Lemma 2.9, if $G$ is a proper subgraph of $H$, $\mu(G)<n-k$, then
$G=K_2\vee(K_{n-k-1}+K_{k-1})$, a contradiction.

If $3\leq l\leq k-1$, we can get $H$ is Hamilton -connected, which
contradicts the assumptions of $H$.

Indeed, let $I$ be the graph induced by $X\cup Y\backslash\{u\}$,
where $u\in Y$. Since $K_{l-1}\vee\ov K_{k-1}\subset I$, and $l\geq
3$, we see that $I$ is 2-connected. Furtherly, if $x$ and $y$ are
distinct nonadjacent vertices of $I$,$$d_I(x)+d_I(y)\geq 2k-2\geq
k+l-1,$$ then $I$ is Hamilton-connected by Lemma 2.1.

Then for any two distinct vertices $x,y$ of $H$, we can get a
Hamilton path of $H$ with $x,y$ as endpoint. So, $H$ is
Hamilton-connected. For example, for any $x,y \in X$. Let $xP_{1}
u_{1}vP_{2}y$ be a Hamilton path of $I$, where $v\in Y$. Let $M$ be
a subgraph of $H$, which is induced by $V(H)\backslash V(I)$. We
note that $M$ is a complete graph, then $M$ is Hamiltonian. So,
there is a Hamilton cycle $C:uP_{3}v_{1}u$ of $M$. Now we delete the
edges $u_{1}v,uv_{1}$, and add the edges $u_{1}u,vv_{1}$, then we
get a path $xP_{1} u_{1}uP_{3} v_{1}vP_{2} y$ be a Hamilton path of
$H$. Similar methods prove the other cases.

If $l=k$, we also can find that $H$ is Hamilton -connected, which
contradicts the assumptions of $H$.  For example, for any $x\in X$,
$y\in Z$. Because every vertex in $Y$ is adjacent to every vertex in
$X$, there is a path $xP_{4} v$, which contains all vertices of
$X\cup Y$, where $v\in Y$. Let $N$ be a subgraph of $H$, which is
induced by $Z\cup \{v\}$. We note that $N$ is a complete graph, then
$N$ is Hamilton-connected. So, there is a Hamilton path $vP_{5} w$
of $N$. Now, we get a path $xP_{4}vP_{5}w$ be a Hamilton path of
$H$. Similar methods prove the other cases.

So, the result follows. \hfill $\blacksquare$

\begin{theorem}
Let $k\geq 1, n\geq2(k+1)^2$, and let $G$ be a graph of order $n$
with minimum degree $\delta(G)\geq k$. If
$$\mu(G)\geq\frac{n^2}{n-1}-\frac{nk}{n-1}-\frac{2}{\sqrt{n-2}},$$
then $G$ is traceable from every vertex, unless
$G=K_1\vee(K_{n-k-1}+K_k)$.
\end{theorem}

{\bf Proof.} Let $H=G\vee K_1$, then $H$ be a graph of order $n+1$,
with minimum degree $\delta(H)\geq k+1$. By Lemma 2.7 and the
assumption. We have
\begin{eqnarray*}
\mu(H)&>&\frac{n-1}{n}\mu(G)+2\frac{\sqrt{n-1}}{n}\\
&\geq&\frac{n-1}{n}(\frac{n^2}{n-1}-\frac{nk}{n-1}-\frac{2}{\sqrt{n-1}})+2\frac{\sqrt{n-1}}{n}\\
&=&(n+1)-(k+1).
\end{eqnarray*}
Then by Theorem 2.10, we get $H$ is Hamilton-connected, unless
$H=K_2\vee(K_{n-k-1}+K_k).$

So, according to the Lemma 2.5, $G$ is traceable from every vertex,
unless $G=K_1\vee(K_{n-k-1}+K_k)$. \hfill $\blacksquare$

Let $ES_n$ be the set of following graphs of even order $n$:

(i) $ K_{\frac{n}{2},\frac{n}{2}}$;

(ii) $G_1\vee G_2$, where $G_1$ is a regular graph of order $n-r$
with degree $\frac{n}{2}-r$, $G_2$ has $r$ vertices, $1\leq
r\leq\frac{n}{2}$.

Let $EW_n$ be the set of following graphs of odd order $n$:

$ G_1\vee G_2$, where $G_1$ is a regular graph of order $n+1-r$ with
degree $\frac{n+1}{2}-r$, $G_2$ has $r-1$ vertices, $1\leq
r\leq\frac{n+1}{2}$.

\begin{theorem}
 Let $G$ be a graph of order $n\geq 2k$, where $k\geq 2$. If $\delta(G)\geq
 k$ and
$$\mu(\ov{G})\leq\sqrt{(k-1)(n-k-1)},$$ then $G$ is Hamilton-connected, unless $G= K_{k-1,n-k-1}\vee
K_2$ or $G= K_{k-1,n-k-1}\vee O_2$ or $G\in ES_n$ and $n=2k$.
\end{theorem}

{\bf Proof.} Let $H=C_{n+1}(G)$. If $H$ is Hamilton-connected, then
so is $G$ by Lemma 2.1. Now we assume that $H$ is not
Hamilton-connected. Note that $H$ is  $(n+1)$-closure of $G$, thus
every two nonadjacent vertices $u$, $v$ of $H$ have degree sum at
most $n$, i.e., $$d_{\ov{H}}(u)+d_{\ov{H}}(v)\geq n-2, ~{\rm for~
any~ edge}~ uv\in E(\ov{H}). \eqno(2.8)$$

Since $d_{G}(u)\geq k$ and $d_{G}(v)\geq k$, we have
$d_{\ov{H}}(u)\leq n-k-1$ and $d_{\ov{H}}(v)\leq n-k-1$. Then
combining (2.8), $k-1 \leq d_{\ov{H}}(u)\leq  n-k-1$, $k-1 \leq
d_{\ov{H}}(v)\leq n-k-1$, this implies that
$$d_{\ov{H}}(u)d_{\ov{H}}(v)\geq
d_{\ov{H}}(u)(n-2-d_{\ov{H}}(u))\geq(k-1)(n-k-1),$$ with equality if
and only if (up to symmetry), $d_{\ov{H}}(u)=k-1$ and
$d_{\ov{H}}(v)=n-k-1$. By Lemma 2.6, Perron-Frobenius theorem, and
the assumption,
$$\sqrt{(k-1)(n-k-1)}\geq\mu(\ov{G})\geq\mu(\ov{H})\geq\min_{uv\in
E(\ov{H})}\sqrt{d_{\ov{H}}(u)d_{\ov{H}}(v)}\geq\sqrt{(k-1)(n-k-1)}.$$

Therefore, $\mu(\ov{G})=\mu(\ov{H})=\sqrt{(k-1)(n-k-1)}$, and
$d_{\ov{H}}(u)+d_{\ov{H}}(v)=n-2$ for any edge $uv\in E(\ov{H})$,
and $d_{\ov{H}}(u)=k-1$, $d_{\ov{H}}(v)=n-k-1$. Note that every
nontrivial component of $\ov{H}$ has a vertex of degree at least
$\frac{n}{2}-1$ and hence of order at least $\frac{n}{2}$. This
implies that $\ov{H}=K_\frac{n}{2}+ K_\frac{n}{2}$ for $n=2k$, or
$\ov{H}$ contains exactly one nontrivial component $F$ which is
either regular or semi-regular, and $\frac{n}{2}\leq|V(F)|\leq n$.

Noting that $\mu(\ov{G})=\mu(\ov{H})$, $\ov{G}\supseteq\ov{H}$, if
$\ov{H}=K_\frac{n}{2}+ K_\frac{n}{2}$ and $n=2k$, then
$\ov{G}=\ov{H}$ by the Perron-Frobenius theorem. So $G=
K_{\frac{n}{2},\frac{n}{2}}\in ES_n$ and $n=2k$, a contradiction.
Therefore we assume that $\ov{H}$ contains exactly one nontrivial
component $F$.

First suppose $F$ is an bipartite semi-regular graph. By the
condition of the degree sum of two adjacent vertices, we have $F$
contains at least $n-2$ vertices. If $F$ contains $n-2$ vertices,
then $\ov{H}=K_{k-1,n-k-1}+ O_2$. Noting that
$\mu(\ov{G})=\mu(\ov{H})$, $\ov{G}\supseteq\ov{H}$, then
$\ov{G}=\ov{H}~{\rm or}~K_{k-1,n-k-1}+ K_2$ by the Perron-Frobenius
theorem. So $G=(K_{k-1}+K_{n-k-1})\vee K_2~{\rm
or}~(K_{k-1}+K_{n-k-1})\vee O_2$, a contradiction. If $F$ contains
$n-1$ vertices. Let $F$ with two partite sets ${X,Y}$, then
$|X|=k-1,|Y|=n-k$ or $|X|=k,|Y|=n-k-1$. Thus according to the edge
number of $F$, we have $(n-k)(k-1)=(k-1)(n-k-1)$ or
$(n-k-1)(k-1)=k(n-k-1)$, a contradiction. If $F$ contains $n$
vertices, let $F$ with two partite sets ${X,Y}$,
 then $|X|=k,|Y|=n-k$ or $|X|=k+1,|Y|=n-k-1$ or $|X|=k-1,|Y|=n-k+1$. If $|X|=k,|Y|=n-k$, according to the
edge number of $F$, we have $(n-k-1)k=(k-1)(n-k)$, $n=2k$, and then
$H=\overline{F}$ is
 Hamilton-connected, a contradiction. If $|X|=k+1,|Y|=n-k-1$ or $|X|=k-1,|Y|=n-k+1$, according to the edge number of $F$, we have
$(n-k-1)(k+1)=(k-1)(n-k-1)$ or $(n-k-1)(k-1)=(k-1)(n-k+1)$, a
contradiction.

Next we assume $F$ is a regular graph. Then for every $v\in V(F)$,
$d_{F}(v)=\frac{n}{2}-1$, and $n=2k$. If $F=\ov{H}$, by a similar
discussion as the above, $\ov{G}=\ov{H}$, and hence $G=H$ is regular
of degree $\frac{n}{2}$. By Lemma 2.8,
$G=K_{\frac{n}{2},\frac{n}{2}}\in ES_{n}$, or $G$ is
Hamilton-connected, a contradiction. Otherwise, $\ov{H}=F\cup O_r$,
where $r=n-|V(F)|$ and $1\leq r\leq \frac{n}{2}$. Noting that
$\mu(\ov{G})=\mu(\ov{H})$, we have $\ov{G}=F\cup F_1$, where $F_1$
is obtained from $O_r$ possibly adding some edges. Hence
$G=\ov{F}\vee\ov F_1\in ES_n$, a contradiction. \hfill
$\blacksquare$

\begin{theorem}
Let $G$ be a graph of order $n\geq 2k+1$, where $k\geq2$. If
$\delta(G)\geq k$ and $$\mu(\ov{G})\leq\sqrt{k(n-k-1)},$$ then $G$
is traceable from every vertex, unless $G= K_{k,n-k-1}\vee K_1$ or
$G\in EW_{n}$ and $n=2k+1$.
\end{theorem}

{\bf Proof.} Let $G'=G\vee K_1$. We note that $|V(G')|=n+1$,
$\mu(\ov{G'})=\mu(\ov{G})\leq\sqrt{k(n-k-1)}$, $\delta(G')\geq k+1$.
By Theorem 2.12, we get $G'$ is Hamilton-connected, unless $G'=
K_{k,n-k-1}\vee K_2$ or $G'=K_{k,n-k-1}\vee O_2$ or $G'\in ES_{n+1}$
and $n=2k+1$. By Lemma 2.5 and the construction of $G'$, we have $G$
is traceable from every vertex, unless $G= K_{k,n-k-1}\vee K_1$ or
$G\in EW_n$ and $n=2k+1$. \hfill $\blacksquare$

\section{(Sigless Laplacian) Spectral radius conditions for a nearly balanced bipartite graph to be traceable}

We note that if a bipartite graph $G=(X,Y;E)$ is traceable, $G$ is a
balanced bipartite graph or a nearly balanced bipartite graph. Li
and Ning \cite{ning1} has presented some (signless Laplacian)
spectral radius conditions for a balanced bipartite graph to be
Hamiltonian. If $G=(X,Y;E)$ be a nearly balanced bipartite graph
with $|X|=|Y|-1$, we can obtained $G'$ from $G$ by adding a vertex
which is adjacent to every vertex in $Y$, then $G'$ be a balanced
partite graph. Note that $G$ is traceable if and only if $G'$ is
Hamiltonian. Inspired by this, in this section, we will study the
conditions for a nearly balanced bipartite graph to be traceable in
terms of spectral radius, signless Laplacian spectral radius of the
graph or its quasi-complement.

Let $G$ be balanced bipartite graph of order $2n$. The {\it
bipartite closure} of $G$, denoted by $cl_{B}(G)$, is the graph
obtained from $G$ by recursively joining pairs of nonadjacent
vertices in different partite sets whose degree sum is at least
$n+1$ until no such pair remains. Note that
$d_{{cl}_B(G)}(u)+d_{{cl}_B(G)}(v)\leq n$ for any pair of
nonadjacent vertices $u$ and $v$ in the distant partite sets of
$cl_B(G)$.

\begin{lemma} {\rm (Bondy and Chv\'{a}tal \cite{bondy})}
A balanced bipartite graph $G$ is Hamiltonian if and only if
$cl_{B}(G)$ is Hamiltonian.
\end{lemma}

Before introducing our results, we need some notations. In order to
facilitate understanding, in this paper, when we mention a bipartite
graph, we always fix its partite sets, e.g., $O_{n,m}$ and $O_{m,n}$
are considered as different bipartite graphs, unless $m=n$.

Let $G_1, G_2$ be two bipartite graphs, with the bipartition $\{X_1,
Y_1\}$ and $\{X_2, Y_2\}$, respectively. We use $G_1\sqcup G_2$ to
denote the graph obtained from $G_1+ G_2$ by adding all possible
edges between $X_1$ and $Y_2$ and all possible edges between $Y_1$
and $X_2$. We define some classes of graphs as follows:
\begin{eqnarray*}
B_n^k&=&O_{k,n-k}\sqcup K_{n-k,k}             (1\leq k\leq n/2),\\
C_n^k&=&O_{k,n-k}\sqcup K_{n-k-1,k}          (1\leq k\leq n/2).
\end{eqnarray*}
Note that $e(B^k_n)=n(n-k)+k^2$, $e(C^{k}_{n})=n(n-k-1)+k^2$,
$\mu(\widehat{B^k_n})=\mu(\widehat{C^k_n})=\mu(K_{k,n-k})=\sqrt{k(n-k)}$,
and $B^k_n$ is not Hamiltonian, $C^{k}_{n}$ is not traceable. By
Perron-Frobenius theorem, $\mu(B^k_n)>\mu(K_{n,n-k})=\sqrt{n(n-k)}$,
$\mu(C^{k}_{n})> \mu(K_{n,n-k-1})=\sqrt{n(n-k-1)}$.

\begin{center}
\vspace{2mm}
\includegraphics[scale=.7]{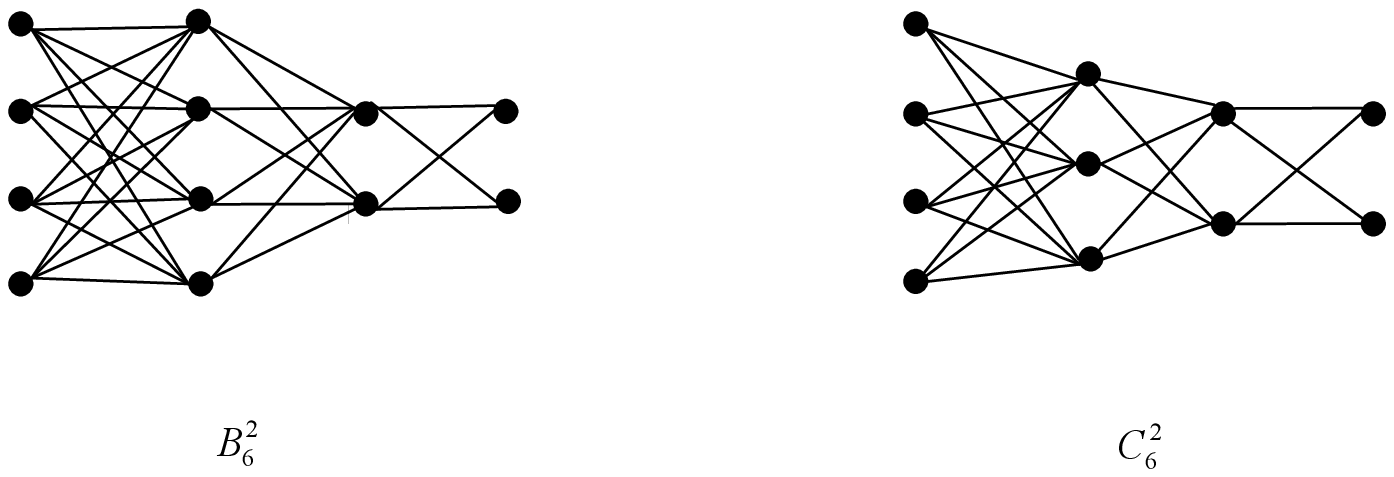}

\vspace{2mm} \small Fig. 3.1. Graphs $B_{6}^{2}$ and $C_{6}^{2}$.
\end{center}

Let $G=(X,Y)$ be a bipartite graph with two part sets $X$, $Y$.
Denote by $\mathscr{B}_{n}^{k}(1\leq k\leq n/2)=\{O_{k,n-k}\sqcup
G(X,Y)$, where $|X|=n-k$, $|Y|=k$\}. Denote by
$\mathscr{C}_{n}^{k}(1\leq k\leq n/2)=\{O_{k,n-k}\sqcup G(X,Y)$,
where $|X|=n-k-1$, $|Y|=k$\}.

\begin{lemma}{\rm (Li and Ning \cite{ning1})}
Let $G$ be a balanced bipartite graph of order $2n$. If
$\delta(G)\geq k\geq 1, n\geq 2k+1$ and $$ e(G)> n(n-k-1)+(k+1)^2,$$
then $G$ is Hamiltonian unless $G\subseteq B^k_n$.
\end{lemma}

\begin{lemma}
Let $G=(X,Y)$ be a nearly balanced bipartite graph of order $2n-1$.
If $\delta(G)\geq k\geq1, n\geq 2k+1$, and
$$e(G)> n(n-k-2)+(k+1)^2,$$ then $G$ is traceable unless $G\subseteq
C^{k}_{n}$.
\end{lemma}

{\bf Proof.} Let $|X|=n-1, |Y|=n$, $G'$ be obtained from $G$ by
adding a vertex which is adjacent to every vertex in $Y$, then $G'$
be a balanced bipartite graph. Note that $G$ is traceable if and
only if $G'$ is Hamiltonian. We have $|V(G')|=2n,
\delta(G')\geq\delta(G)\geq k\geq1, n\geq 2k+1,$ and
$$e(G')=e(G)+n>n(n-k-2)+(k+1)^2+n=n(n-k-1)+(k+1)^2.$$ By lemma 3.2, $G'$ is
Hamiltonian unless $G'\subseteq B^k_n$. Thus $G$ is traceable unless
$G\subseteq C^{k}_{n}$. \hfill $\blacksquare$

\begin{lemma} \rm{(Bhattacharya, Friedland and
Peled \cite{Bbhattacharya})} Let $G$ be a bipartite graph. Then
$$\mu(G)\leq\sqrt{e(G)}.$$
\end{lemma}

\begin{lemma} \rm{(Ferrara, Jacobson and Powell \cite{Ferrara})}
Let $G$ be a non-Hamiltonian balanced bipartite graph of order $2n$.
If $d(u)+d(v)\geq n$ for every two nonadjacent vertices $u,v$ in
distinct partite sets, then either $G\in
\bigcup_{k=1}^{n/2}\mathscr{B}_{n}^{k}$, or $G=\Gamma_{1}$ or
$\Gamma_{2}$ for $n=4$.
\end{lemma}

\begin{center}
\vspace{2mm}
\includegraphics[scale=.7]{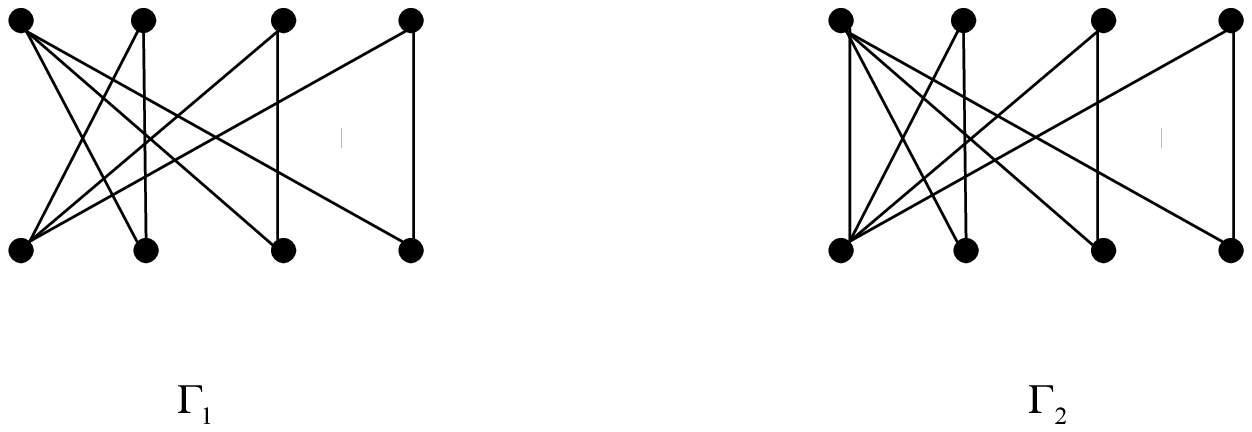}

\vspace{2mm} \small Fig. 3.2. Graphs $\Gamma_{1}$ and $\Gamma_{2}$.
\end{center}

\begin{lemma}{\rm(Feng and Yu \cite{feng}, Yu and Fan \cite{yu})}
Let G be a graph with non-empty edge set. Then
$$q(G)\leq \max\left\{d(u)+\frac{\sum\limits_{v\in N(u)}d(v)}{d(u)}: u\in V(G)\right\}.$$
\end{lemma}

\begin{lemma}
 Let $G$ be a bipartite graph with two partite sets $X,~Y$, and $\max
 \{|X|,~|Y|\}=n$. Then $$q(G)\leq \frac{e(G)}{n}+n.$$
\end{lemma}

{\bf Proof.} If $G$ is an edgeless graph, then $q(G)=0$, and the
result is trivially true. Now assume $G$ contains at lest one edge.
Let $x\in V(G)$, and
$$d(x)+\frac{\sum\limits_{v\in N(x)}d(v)}{d(x)}=\max\left\{d(u)+\frac{\sum\limits_{v\in N(u)}d(v)}{d(u)}: u\in V(G)\right\}.$$
By Lemma 3.6 and for every $v\in V(G),d_G(v)\leq\max
 \{|X|,~|Y|\}=n$, we get
$$\begin{array}{ll} \frac{e(G)}{n}+n-q(G) &\geq\left(\frac{\sum\limits_{v\in
N(x)}d(v)}{n}+n\right)-\left(d(x)+\frac{\sum\limits_{v\in
N(x)}d(v)}{d(x)}\right)\\&=(n-d(x))\left(1-\frac{\sum\limits_{v\in
N(x)}d(v)}{nd(x)}\right)\\
&\geq 0.
\end{array}.$$

The result follows. \hfill $\blacksquare$

\begin{lemma} Let $k\geq 1, n\geq\frac{k^3}{2}+k+2$. If $G$ is a subgraph of $C_n^k$, $\delta(G)\geq
k$. Then $\mu(G)<\sqrt{n(n-k-1)}$, unless $G=C_n^k$.
\end{lemma}

{\bf Proof.} The proof is similar to Lemma 2.9. Set for short $\mu
:=\mu(G)$, and let $\textbf{x}=(x_{v_1},\ldots,x_{v_{2n-1}})^{T}$ be
a unit Perron vector of $G$. By (2.3), we have
$$\mu=\textbf{x}^TA(G)\textbf{x}.$$

Assume that $G$ is a proper subgraph of $C_n^k$. By Perron-Frobenius
theorem, we may assume that $G$ is obtained by omitting just one
edge $uv$ of $C_n^k$.

Write $X$ for the set of vertices of $C_n^k$ of degree $k$, let $Y$
be the set of vertices of $C_n^k$ of degree $n$, let $Z$ for the set
of vertices of $C_n^k$ of degree $n-k-1$, let $H$ be the set of the
remaining $k$ vertices of $C_n^k$ of degree $n-1$.

Since $\delta(G)\geq k$, we can see that $G$ must contain all the
edges between $X$ and $H$. Therefore $\{u,v\}\subset Y\cup H$ or
$\{u,v\}\subset Y\cup Z$, with two possible cases: $(a)~ u\in Y,v\in
H$; $(b)~ u\in Y, v\in Z$. We shall show that case $(b)$ yields a
graph of no smaller spectral radius than case $(a)$.

Indeed, by (2.2), we have $x_i=x_j$ for any $i,j\in X$; likewise
$x_i=x_j$ for any $i,j\in Y\backslash\{u\}$, for any $i,j\in
Z\backslash\{v\}$ and for any $i,j\in H\backslash\{v\}$. Thus, let
\begin{eqnarray*}
x &:=& x_i, i\in X,\\
y &:=& y_i, i\in Y\backslash\{u\},\\
z &:=& z_i, i\in Z\backslash\{v\},\\
h &:=& h_i, i\in H\backslash\{v\}.
\end{eqnarray*}

Suppose that case $(a)$ holds, that is, $u\in Y, v\in H$. Choose a
vertex $w\in Z$, remove the edge $uw$, and add the edge $uv$. Then
the obtained graph $G'$ is covered by case $(b)$.

If $x_w\leq x_v$, we have
$$\textbf{x}^TA(G')\textbf{x}-\textbf{x}^TA(G)\textbf{x}=2x_u(x_v-x_w)\geq0;$$

If $x_w>x_v$, swap the entries $x_v$ and $x_w$, write $\textbf{x}'$
for the resulting vector. We note that $\textbf{x}'$ is also a unit
vector, and have that
$$\textbf{x}'^TA(G')\textbf{x}'-\textbf{x}^TA(G)\textbf{x}=2(x_w-x_v)\sum_{i\in X}x_i\geq0.$$ Then
by (2.3), $\mu(G')\geq\mu(G)$, as claimed.

Therefore, we may assume that $u\in Y, v\in Z$, and set $t:=x_u,
s:=x_v$, note that the $2n-1$ eigenvalue-equations of $G$ are
reduced to six equations involving just the unknowns $x,y,z,h,t,$
and $s$:
\begin{eqnarray*}
\mu x&=&kh,\\
\mu y&=&(n-k-1)z+kh+s,\\
\mu z&=&(n-k-2)y+t,\\
\mu h&=&kx+(n-k-2)y+t,\\
\mu t&=&(n-k-1)z+kh,\\
\mu s&=&(n-k-2)y.
\end{eqnarray*}
We find that
\begin{eqnarray*}
x&=&\frac{k}{\mu}h,\\
t&=&\frac{(n-k-1)(\mu^2-k^2)+k\mu^2}{\mu^3}h,\\
s&=&\frac{(\mu^2-(n-k-1))(\mu^2-k^2)-k\mu^2}{\mu^4}h.
\end{eqnarray*}

Furtherly, note that if we remove all edges between $X$ and $H$, and
add the edge $uv$ to $G$, we obtain the graph $K_{n,n-k-1}+\ov K_k$.
Letting $\textbf{x}''$ be the restriction of $\textbf{x}$ to
$K_{n,n-k-1}$, we find that
$$\textbf{x}''^TA(K_{n,n-k-1})\textbf{x}''=\textbf{x}^TA(G)\textbf{x}+2st-2k^2xh=\mu+2st-2k^2xh.$$ But since
$\|\textbf{x}''\|^2=1-kx^2$, we see that
$$\mu+2st-2k^2xh=\textbf{x}''^TA(K_{n,n-k-1})\textbf{x}''\leq\mu(K_{n,n-k-1})\|\textbf{x}''\|^2=\sqrt{n(n-k-1)}(1-kx^2).$$

Assume for a contradiction that $\mu\geq\sqrt{n(n-k-1)}$. This
assumption together with above inequality, yields
$$\mu+2st-2k^2xh\leq\mu(1-kx^2),$$ and therefore
$$2st-2k^2xh\leq-kx^2\mu.$$
Now, first combining above equality about $x$, then combining above
equalities about $t$ and $s$, we have
$$k^3\geq\frac{2\mu}{h^2}st\geq\frac{2(n-k-1)(\mu^2-(n-k-1))(\mu^2-k^2)^2}{\mu^6}-\frac{2k(n-k-1)(\mu^2-k^2)}{\mu^4}.$$
Applying Bernoulli's inequality to the right side, we get
\begin{eqnarray*}
k^3&\geq&2(n-k-1)(\frac{\mu^2-(n-k-1)}{\mu^2})(\frac{\mu+k}{\mu})^2(\frac{\mu-k}{\mu})^2-\frac{2k(n-k-1)}{\mu^2}(1-\frac{k^2}{\mu^2})\\
&=&2(n-k-1)(1-\frac{n-k-1}{\mu^2})(1+\frac{k}{\mu})^2(1-\frac{k}{\mu})^2-\frac{2k(n-k-1)}{\mu^2}+\frac{2k^3(n-k-1)}{\mu^4}\\
&>&2(n-k-1)(1-\frac{n-k-1}{\mu^2})-\frac{2k(n-k-1)}{\mu^2}.
\end{eqnarray*}
Using the inequality $\mu\geq\sqrt{n(n-k-1)}$, we easily find that
$$k^3>2(n-k-1)-\frac{2(n-1)}{n}>2(n-k-2),$$  and then $n<\frac{k^3}{2}+k+2$, a
contradiction. \hfill $\blacksquare$

\begin{theorem}
Let $G$ be a nearly balanced bipartite graph of order $2n-1
~(n\geq\max\{\frac{k^3}{2}+k+2,(k+1)^2\})$, where $k\geq1$. If
$\delta(G)\geq k$ and
$$\mu(G)>\sqrt{n(n-k-1)},$$ then $G$ is traceable, unless $G=C_n^k$.
\end{theorem}

{\bf Proof.} By the assumption and Lemma 3.4,
$$\sqrt{n(n-k-1)}<\mu(G)\leq\sqrt{e(G)}.$$ Thus, we obtain
$$e(G)>n(n-k-1)\geq n(n-k-2)+(k+1)^2,$$ when $n\geq\max\{\frac{k^3}{2}+k+2,(k+1)^2\}>2k+1$, by Lemma
3.3, $G$ is traceable or $G\subseteq C^{k}_{n}$. But if $G\subseteq
C^{k}_{n}$, then $\mu(G)<\sqrt{n(n-k-1)}$, unless $G=C_n^k$ by Lemma
3.8, a contradiction. \hfill $\blacksquare$

\begin{theorem}
Let $G=(X,Y)$ be a nearly balanced bipartite graph of order $2n-1
~(n\geq2k)$, where $k\geq1$. If $\delta(G)\geq k$, and
$$\mu(\widehat{G})\leq\sqrt{k(n-k)},$$ then $G$ is traceable,
unless $G\in \bigcup_{k=1}^{n/2}\mathscr{C}^{k}_{n}$ or
$\Gamma_{2}-v$, where $d_{\Gamma_{2}}(v)=4$.
\end{theorem}

{\bf Proof.} Let $|X|=n-1, |Y|=n$, $G'$ be obtained from $G$ by
adding a vertex which is adjacent to every vertex in $Y$, then $G'$
be a balanced partite graph. Note that $G$ is traceable if and only
if $G'$ is Hamiltonian. Let $H=cl_{B}(G')$. If $H$ is Hamiltonian,
then so is $G'$ by Lemma 3.1. Now we assume that $H$ is not
Hamiltonian. Note that $H$ is bipartite closure of $G$, thus every
two nonadjacent vertices $u$, $v$ in distant part sets of $H$ have
degree sum at most $n$, i.e.,
$$d_{\widehat{H}}(u)+d_{\widehat{H}}(v)=n-d_{H}(u)+n-d_{H}(v)\geq
n, \mbox{~for any edge~} uv\in E(\widehat{H}). \eqno (3.1)$$

This implies that $\widehat{H}$ contains only one component or
$\widehat{H}=K_{s,n-s}+ K_{t,n-t},~ s,t\geq 1$. If
$\widehat{H}=K_{s,n-s}+ K_{t,n-t},~ s,t\geq 1$, it contradicts the
structure of $\widehat{H}$ (It must contains an isolated vertex).
So, $\widehat{H}$ contains only one component.

Since $\delta(H)\geq\delta(G^{'})\geq\delta(G)\geq k$, we can see
that $d_{\widehat{H}}(u)\leq n-k$ and $d_{\widehat{H}}(v)\leq n-k$.
Thus by (3.1), we have $k\leq d_{\widehat{H}}(u)\leq n-k$, $k\leq
d_{\widehat{H}}(v)\leq n-k$, this implies that
$$d_{\widehat{H}}(u)d_{\widehat{H}}(v)\geq
d_{\widehat{H}}(u)(n-d_{\widehat{H}}(u))\geq k(n-k),$$ with equality
if and only if (up to symmetry) $d_{\widehat{H}}(u)=k$,
$d_{\widehat{H}}(v)=n-k$. By Lemma 2.6,
$$\sqrt{k(n-k)}\geq\mu(\widehat{G})=\mu(\widehat{G^{'}})\geq
\mu(\widehat{H})\geq \min_{uv\in
E(\widehat{H})}\sqrt{d_{\widehat{H}}(u)d_{\widehat{H}}(v)}\geq
\sqrt{k(n-k)},$$ this implies that $\mu(\widehat{H})=\sqrt{k(n-k)}$
and there is an edge $uv\in E(\widehat{H})$ such that
$d_{\widehat{H}}(u)=k$, $d_{\widehat{H}}(v)=n-k$. Let $F$ be the
component of $\widehat{H}$ which contains $uv$. By Lemma 2.6, $F$ is
an bipartite semi-regular graph, with partite sets $X'\subseteq X$,
and $Y'\subseteq Y$, and for any vertex $x\in X'$, $d_F(x)=k$, and
any vertex $y\in Y'$, $d_F(y)=n-k$. Then $d_H(u)+d_H(v)=n$ for every
two nonadjacent vertices $u,v$ in distinct partite sets of $H$. By
Lemma 3.5, $H\in \bigcup_{k=1}^{n/2}\mathscr{B}_{n}^{k}$ or
$H=\Gamma_{1}$ or $\Gamma_{2}$ for $n=4$ and $k=2$, then
$G^{'}\subseteq B_{n}^{k}$ ($1\leq k\leq n/2$) or $G^{'}\subseteq
\Gamma_{1}$ or $\Gamma_{2}$ for $n=4$ and $k=2$. By Perron-Frobenius
theorem, every (spanning) subgraph of $\Gamma_{1}$, $\Gamma_{2}$ or
$B_{n}^{k}$, $1\leq k\leq n/2$, if it is not $\Gamma_{1}$ or
$\Gamma_{2}$ or a graph in $\mathscr{B}_{n}^{k}$, $1\leq k\leq n/2$,
then has the quasi-complement with spectral radius greater than
$\sqrt{k(n-k)}$. Thus $G^{'}\in
\bigcup_{k=1}^{n/2}\mathscr{B}_{n}^{k}$ or $\Gamma_{1}$ or
$\Gamma_{2}$ for $n=4$ and $k=2$. By the construction of $G^{'}$, we
get $G\in \bigcup_{k=1}^{n/2}\mathscr{C}^{k}_{n}$ or $\Gamma_{2}-v$,
where $d_{\Gamma_{2}}(v)=4$, a contradiction.

\begin{theorem}
Let $G$ be a nearly balanced bipartite graph of order $2n-1
(n\geq(k+1)^2)$, where $k\geq1$. If $\delta(G)\geq k$ and
$$q(G)>\frac{n(2n-k-2)+(k+1)^{2}}{n},$$
then $G$ is traceable, unless $G\subseteq C_{n}^{k}$.
\end{theorem}

{\bf Proof.} By the assumption and Lemma 3.7,
$$\frac{n(2n-k-2)+(k+1)^{2}}{n}<q(G)\leq\frac{e(G)}{n}+n.$$ Thus, we obtain
$$e(G)>n(n-k-1)\geq n(n-k-2)+(k+1)^2,$$ when $n\geq(k+1)^2$, by Lemma
3.3, $G$ is traceable or $G\subseteq C^{k}_{n}$.  \hfill
$\blacksquare$

{\bf Remark:} In Theorem 3.11, we can't change $G\subseteq
C_{n}^{k}$ to $G= C_{n}^{k}$ like Theorem 3.9. In fact, we can find
a subgraph $G\subset C_{n}^{k}$, which satisfies the conditions of
Theorem 3.11, such tat
$q(G)>2n-k-1\geq\frac{n(2n-k-2)+(k+1)^{2}}{n}$.

{\bf Proof.} Assume that $G$ is a proper subgraph of $C_n^k$,
$n\geq(k+1)^2, k\geq1, \delta(G)\geq k$, and has the maximum
signless Laplacian spectral. By Perron-Frobenius theorem, $G$ is
obtained by omitting just one edge $uv$ of $C_n^k$. Set for short $q
:=q(G)$, and let $\textbf{x}=(x_{v_1},\ldots,x_{v_{2n-1}})^{T}$ be a
positive unit eigenvector to $q$. We have
$$q=\textbf{x}^TQ(G)\textbf{x}=2\sum\limits_{uv\in
E(G)}(x_u+x_v)^2,\eqno(3.2)$$ and $$(q-d_G(v))x_v=\sum\limits_{u\in
N_G(v)}x_u, \eqno(3.3)$$ for each vertex $v\in V(G)$. Equation (3.3)
is called the signless Laplacian eigenvalue-equation for the graph
$G$. In addition, for an arbitrary unit vector $\textbf{x}\in
R^n$,$$q\geq \textbf{x}^TQ(G)\textbf{x}, \eqno(3.4)$$ with equality
holds if and only if $\textbf{x}$ is an eigenvector of $Q(G)$
according to $q$.

Write $X$ for the set of vertices of $C_n^k$ of degree $k$, let $Y$
be the set of vertices of $C_n^k$ of degree $n$, let $Z$ for the set
of vertices of $C_n^k$ of degree $n-k-1$, let $H$ be the set of the
remaining $k$ vertices of $C_n^k$ of degree $n-1$.

Since $\delta(G)\geq k$, we see that $G$ must contain all the edges
between $X$ and $H$. Therefore $\{u,v\}\subset Y\cup H$ or
$\{u,v\}\subset Y\cup Z$, with two possible cases: $(a)~ u\in Y,v\in
H$; $(b) ~u\in Y, v\in Z$. We shall show that case $(b)$ yields a
graph of no smaller signless Laplacian spectral radius than case
$(a)$.

Indeed, by (3.3), we have $x_i=x_j$ for any $i,j\in X$; likewise
$x_i=x_j$ for any $i,j\in Y\backslash\{u\}$, for any $i,j\in
Z\backslash\{v\}$ and for any $i,j\in H\backslash\{v\}$. Thus, let
\begin{eqnarray*}
x &:=& x_i, i\in X,\\
y &:=& y_i, i\in Y\backslash\{u\},\\
z &:=& z_i, i\in Z\backslash\{v\},\\
h &:=& h_i, i\in H\backslash\{v\}.
\end{eqnarray*}

Suppose that case $(a)$ holds, that is, $u\in Y, v\in H$. Choose a
vertex $w\in Z$, remove the edge $uw$, and add the edge $uv$. Then
the obtained graph $G'$ is covered by case $(b)$.

If $x_w\leq x_v$, we have
$$\textbf{x}^TQ(G')\textbf{x}-\textbf{x}^TQ(G)\textbf{x}=2(x_v+x_u)^2-2(x_w+x_u)^2\geq0;$$

If $x_w>x_v$, swap the entries $x_v$ and $x_w$, write $\textbf{x}'$
for the resulting vector. We note that $\textbf{x}'$ is a unit
vector, and have that
$$\textbf{x}'^TQ(G')\textbf{x}'-\textbf{x}^TQ(G)\textbf{x}=2(x_w+\sum_{i\in X}x_i)^2-2(x_v+\sum_{i\in X}x_i)^2\geq0.$$
Then by (3.4), $q(G')\geq q(G)$.

Therefore, $G$ is obtained by omitting just one edge $uv$ of
$C_{n}^{k}$, where $u\in Y, v\in Z$. Now set $t:=x_u, s:=x_v$, note
that the $2n-1$ signless Laplacian eigenvalue-equations of $G$  are
reduced to six equations involving just the unknowns $x,y,z,h,t,$
and $s$. By Equation (3.3), we have
\begin{eqnarray*}
(q-k)x &=& kh,\\
(q-n)y &=& (n-k-1)z+kh+s,\\
(q-(n-k-1))z &=& (n-k-2)y+t,\\
(q-(n-1))h &=& kx+(n-k-2)y+t,\\
(q-(n-1))t &=& (n-k-1)z+kh,\\
(q-(n-k-2))s &=& (n-k-2)y. \end{eqnarray*} Transform the above
equations into a matrix equation $(B-qI)\textbf{x}=0$, where
$\textbf{x}=(x,y,z,h,t,s)^T$
\begin{equation}
B=\left(
\begin{array}{cccccc}
k&0&0&k&0&0\\
0&n&n-k-1&k&0&1\\
0&n-k-2&n-k-1&0&1&0\\
k&n-k-2&0&n-1&1&0\\
0&0&n-k-1&k&n-1&0\\
0&n-k-2&0&0&0&n-k-2\\
\end{array}
\right)\nonumber.
\end{equation}

Let
$$\begin{aligned}
f(x)&=\det(B-xI)\\&=-x(n-1-x)\left(x^4+(-4n+k+4)x^3+(-nk+6+5n^2-2k^2-11n+k)x^2\right.\\
&\left.+(7n^2+5nk+2+6nk^2-2n^2k-7n-2n^3-6k^2-2k^3-3k)x\right.\\
&\left.+2nk^3-k^3-2k-3k^2+8nk^2+7nk+2n^3k-4n^2k^2-7n^2k\right).
\end{aligned}$$

Thus, $q$ is the largest root of $f(x)=0$, and when $x>q$, $f(x)$ is
is monotonically increasing.

But when $k\geq1, n\geq(k+1)^2,$ we have
$$f(2n-k-1)=(k+1-2n)(k-n)((4-k^2)n^2+(-6k+k^3-4)n+1+3k^2+3k+k^3)<0,$$
which implies that $q>2n-k-1$, the result follows.

\end{document}